\documentclass[12pt,reqno]{amsart}
\usepackage{soul}
\usepackage{multicol,amsmath,graphics, mathtools, cancel,soul,color,bbm,amsfonts,dsfont,tikz,mathrsfs,amssymb,bm,bbold,comment}
\usepackage{algorithm}
\usepackage{algorithmicx}
\usepackage{algpseudocode}

\usepackage[foot]{amsaddr}

\usepackage{xcolor}
\usepackage[labelformat=parens,labelsep=space]{subcaption}

\usepackage{graphicx}

\usepackage{stmaryrd}
\usepackage[left=1in, right=1in, top=1.1in,bottom=1.1in]{geometry}
\usetikzlibrary{matrix,patterns,positioning}
\numberwithin{equation}{section}

\usepackage[normalem]{ulem} % JHRG: Crossedout text
\newcommand{\mathsout}[1]{\ifmmode\text{\sout{\ensuremath{#1}}}\else\sout{#1}\fi}

\usepackage{comment}
\usepackage{soul}

\usepackage[hidelinks]{hyperref}

\newcommand{\E}{\mathbb{E}}

\newcommand{\N}{\mathbb{N}}
\newcommand{\Pb}{\mathbb{P}}

\newcommand{\R}{\mathbb{R}}

\newcommand{\floor}[1]{\left\lfloor#1\right\rfloor}
\newcommand{\vertiii}[1]{{\left\vert\kern-0.25ex\left\vert\kern-0.25ex\left\vert #1
    \right\vert\kern-0.25ex\right\vert\kern-0.25ex\right\vert}}

\def\R{\mathbb{R}}

\def\dh2l{\mathbf{d}_{\mathbb{H}_{2\ell}}}
\def\d2{\mathbf{d}_2}

\newcommand{\Indi}[1]{\mathbbm{1}_{{#1}}}

\newtheorem{theorem}{Theorem}[section]

\newtheorem{proposition}[theorem]{Proposition}

\theoremstyle{remark}

\theoremstyle{definition}

 \numberwithin{dummy}{section}

\def\1{\mathbbm{1}}

\makeatletter
\@namedef{subjclassname@1991}{2020 Mathematics Subject Classification}
\makeatother

\begin{document}

\title{Additive functionals of Harmonic samples: the conditioned Dickman regime}
%\author{Arturo Jaramillo \and Isa\'ias BaÃ±ales \and Joshu\'e Hel\'i Ricalde-Guerrero }
%\date{\today}

\author[Victor Bernal Ramirez]{Victor Bernal Ramirez}
\address{$\dag$) University of California, Santa Barbara}
\email{$\dag$) vbernalramirez@ucsb.edu}

\author[Arturo Jaramillo]{Arturo Jaramillo}
\address{$\ast$) Centro de Investigaci\'on en Matem\'aticas A.C., Department of Probability and Statistics, M\'exico} 
\email{$\ast$)jagil@cimat.mx}

\begin{abstract}
We study the distributional behavior of additive arithmetic functions evaluated at integers drawn from the harmonic distribution. Our main result shows that, for a broad family of completely additive functions, their evaluations at harmonic samples, suitably normalized, converge in law to conditioned Dickman-type Poisson integrals. This behavior contrasts with the Gaussian limits arising in the classical Erd\"os-Kac theorem under uniform sampling. Our approach combines the probabilistic representation of harmonic samples via independent geometric variables, analytic inputs such as Mertens' approximation, and a Poissonization procedure.\\


\end{abstract}

\maketitle

%\tableofcontents

%%%%%%%%%%%%%%%%%%%%%%%%%%%%%%%%%%%%%%%%%%%%%%%%%%%%%%%%%%%%%%%%%%%%%%%%%

%%%%%%%%%%%%%%%%%%%%%%%%%%%%%%%%%%%%%%%%%%%%%%%%%%%%
%   Introduccion
%%%%%%%%%%%%%%%%%%%%%%%%%%%%%%%%%%%%%%%%%%%%%%%%%%%%

\section{Introduction}\label{sec:intro}
\noindent Let $\varphi:\mathbb{N}\to\mathbb{R}$ be a {completely additive} arithmetic function. That is to say, a function satisfying
\[
\varphi(ab)=\varphi(a)+\varphi(b)
\]
for all positive integers $a$ and $b$. In particular,
\[
\varphi(p^\ell)=\ell\,\varphi(p)
\]
for every prime $p$ and every $\ell\geq 0$. This paper is devoted to the study of the behavior of additive functions evaluated at random integers drawn from a {harmonic distribution}. More precisely, let $H_n$ be a random variable taking values in $[n]:=\{1,\ldots,n\}$, with probabilities
\[
\Pb[H_n=k] = \frac{1}{kL_n}\Indi{[n]}(k), \]
where $L_n=\sum_{k=1}^n k^{-1}$. Our main contribution consists of exhibiting a family of completely additive functions $\varphi$ for which suitable normalizations of $\varphi(H_n)$ converge weakly to conditional Dickman-type Poisson integral limits. This type of limiting behavior connects with a classical and well-developed line of work on the probabilistic structure of additive functions and their limiting distributions. Before stating our main theorems precisely, it is therefore natural to situate our contribution within the broader context of probabilistic number theory.\\

\noindent We begin by recalling the classical framework from which this line of results originates. Our starting point is the celebrated Erd\"os-Kac theorem \cite{ErdosKac1940}, which describes the statistical behavior of the prime factor counting function under uniform sampling from $[n]$. More precisely, let $\mathcal{P}$ denote the set of prime numbers and let $\mathcal{P}_n$ be defined as those elements in $\mathcal{P}$ bounded by $n$. The Erd\"os-Kac theorem studies the additive function
\[
\omega(k) = \big|\{\,p\in\mathcal{P} : p \text{ divides } k\,\}\big|.
\]
By means of a truncation argument combined with the method of moments, it was proved that
\[
\frac{\omega(J_n)-\log\log n}{\sqrt{\log\log n}}
\;\;\stackrel{Law}{\longrightarrow}\;\;
N,
\]
where $J_n$ is uniformly distributed on $[n]$ and $N$ is a standard normal random variable.   This theorem is a cornerstone of probabilistic number theory. Beyond its intrinsic statement, it initiated the formal use of probabilistic models to study arithmetic functions. In particular, it introduced the systematic use of stochastic independence arguments in analytic number theory, a perspective that was subsequently developed in depth by Kubilius and later by Elliott (see~\cite{Elliott1979, Elliott1980, Kubilius1964}). Building on this framework, much of the subsequent research was focused on analyzing 
the quantitative content of the works by Kubilius, Elliott, Erd\"os, and Kac. In particular, the convergence rates in the Erd\"os-Kac theorem have been refined in a long line of work, including contributions by R\'enyi and Tur\'an~\cite{RenyiTuran1958}, 
Kowalski and Nikeghbali~\cite{KowalskiNikeghbali2014}, and more recently Chen, Jaramillo, and Yang~\cite{ChenJaramilloYang2021}. On a different front, a parallel line of work, emerging soon after the introduction of probabilistic number theory, focused on extending existing limit theorems to functional settings, where additive functions are viewed as evolving stochastic processes.\\

\noindent The present work remains closer in spirit to classical convergence results, 
without addressing functional or distance-based formulations,  although it draws substantial motivation from both fronts. Among these, we will emphasize the contribution of Billingsley, presented in his Wald Lectures and later published in 1974~\cite{Billingsley1974}.  In that work, Billingsley compared the classical model of uniform sampling of integers from $[n]$ with the Kubilius model (see \cite{Kubilius1964}), in which the event ``$p$ divides $n$'' is replaced by an independent Bernoulli random variable of mean $1/p$. This reformulation embeds additive functions into a probabilistic setting where classical limit theorems apply directly. Within this framework, Billingsley showed that, once centered and normalized, additive functions converge to Gaussian laws under Lindeberg-type conditions, but may also give rise to more general infinitely divisible limits when prime contributions have heavier tails. He further established functional central limit theorems: in the Gaussian regime the limiting process is Brownian motion, while in the general case it is a L\'evy process.\\

\noindent While the uniform sample and Kubilius model form the natural starting point for probabilistic number theory, the present work turns from uniform sampling to the closely related harmonic distribution $H_n$, recently examined in the study of additive functions under alternative sampling frameworks. This model exhibits a clean probabilistic structure, as it admits a stochastic representation in terms of independent geometric random variables  (see~\cite{JaramilloYang2023} and Section~\ref{sec:prelim}), which facilitates the use of standard tools from distributional limit theorems, such as the Stein method.  Its effectiveness is best illustrated in problems such as the description of the law of the largest prime factor of $H_n$, where the perspective uncovers deep connections with the Dickman distribution (see~\cite{Dickman1930} for the original work in the uniform setting)  and yields quantitative approximations with nearly optimal rates, both for the uniform sampling model~\cite{BhattacharjeeSchulte2022} and for the harmonic model~\cite{JaramilloYang2023}.\\

\noindent Returning to the central problem: our goal is to describe the limiting law of suitable normalizations of $\varphi(H_n)$. Our approach addresses this by combining the probabilistic representation of harmonic samples introduced in~\cite{JaramilloYang2023} with classical analytic inputs, most notably Mertens' approximation. In this framework, the harmonic variable $H_n$ admits a representation in terms of independent geometric variables, which makes it natural to study additive functions through their characteristic functions. To analyze these effectively, we introduce a Poissonization procedure that expresses the relevant quantities as integrals against a Poisson point process, thereby ensuring tractability in the study of the associated Fourier transforms. Finally, Mertens-type estimates identify the limiting exponent, which in turn yields the Dickman-type Poisson integral appearing in the limit.\\

\noindent The rest of the paper is organized as follows. Section~\ref{sec:prelim} collects the probabilistic and analytic ingredients of our approach, including a representation of the harmonic distribution via independent geometric variables, a Poissonization procedure that embeds our problem into a stochastic representation of prime multiplicities, classical Mertens-type estimates for sums over primes, and a brief review of the Dickman distribution together with its realization as a Poisson integral. Section~\ref{sec:main} presents the statement of our main results, introducing a new class of limit theorems for additive functions under harmonic sampling. Sections~4 and~5 contain the proofs of the main results.

\section{Preliminaries}\label{sec:prelim}
\noindent In this section we collect the main probabilistic and analytic ingredients underlying our approach. We begin by describing the stochastic representation of the harmonic distribution in terms of independent geometric variables, followed by a Poissonization procedure that enables the analysis of the geometric random variables via characteristic functions. We then review a set of estimates of arithmetic averages over $\mathcal{P}$, derived from Mertens' formulas, and conclude with a brief account of the Dickman distribution and the framework of infinitely divisible laws, which together provide the natural limiting setting for additive arithmetic functions.

\subsection{Prime multiplicities of harmonic samples}
A natural way to describe the arithmetic structure of an integer is through its prime factorization, where every number can be uniquely written as a product of prime powers. For each prime $p$, we record the multiplicity of $p$ in this decomposition through the $p$-adic valuation $\alpha_p:\mathbb{N}\to\mathbb{N}$, defined by
\[
\alpha_p(k) = \max\{\,m \ge 0 : p^{m} \mid k\,\}.
\]
\noindent When integers are sampled uniformly from $[n]$, the $p$-adic valuations exhibit a remarkably regular probabilistic behavior: asymptotically, each $\alpha_p(k)$ follows a geometric distribution, and any finite collection of such multiplicities becomes asymptotically independent. To make this connection precise, let $\varepsilon_p$, with $p$ in $\mathcal P_n$,  be a family of independent random variables with
\[
\Pb[\varepsilon_p = m] = (1 - 1/p)\,p^{-m},
\]
for $m$ a non-negative integer. Then, for each fixed prime $p$, the distribution of $\alpha_p(k)$ under uniform sampling on $[n]$ converges in law to that of $\varepsilon_p$. This  asymptotic picture observed under uniform sampling admits an exact analogue, up to a suitable adjustment, in the harmonic model. In this setting, each prime multiplicity follows an exact geometric law, with dependence arising  from a single global conditioning event. This correspondence, formalized in~\cite[pp.~18-20]{JaramilloYang2023}, establishes that, for any completely additive function $\varphi$,
\begin{equation} \label{eq:HN_decomposition}
\mathcal{L}(\varphi(H_{n})) =
\mathcal{L}\Big(\sum_{p \in \mathcal{P}_{n}} \varphi(p^{\varepsilon_{p}})\,\Big|\,A_{n}\Big),
\end{equation}
where
\[
A_{n} := \Big\{\, \prod_{p \in \mathcal{P}_{n}} p^{\varepsilon_{p}} \le n \,\Big\}.
\]
\noindent Classical estimates derived from Mertens' theorem (see~\cite{ChenJaramilloYang2021}) imply that
\[
\lim_{n\to\infty}\mathbb{P}[A_{n}] = e^{-\gamma},
\]
where $\gamma$ denotes the Euler-Mascheroni constant. Moreover, for all $n \ge 21$, one has $\mathbb{P}[A_{n}] \ge 0.5$.  As a consequence, the conditioning event $A_n$ introduces only a mild distortion of the distribution, and the effective dependence of prime multiplicities remains tractable.\\

\subsection{Poissonization and linear approximation}\label{eq:epsilonpdef}
Building on the geometric decomposition of $H_n$, the next step is to develop a representation that allows one to handle the dependence structure introduced by the conditioning event in~\eqref{eq:HN_decomposition}. To achieve this, we introduce a Poissonization procedure that encodes the prime exponents as points of a Poisson process, yielding a  tractable  method for studying the associated additive functionals. Concretely, we consider the  discrete state space 
\[
\mathcal{Y} := \{(p,k) : p \in \mathcal{P}, k \in \mathbb{N}\},
\]
and a Poisson point process $N$ on $\mathcal{Y}$ with intensity $\nu(p,k) = 1/(k p^k)$, so that the geometric variable is represented as
\begin{align}\label{eq:epsilondef}
\varepsilon_p \;\stackrel{\text{law}}{=}\; \sum_{k \ge 1} k\, N(p,k).	
\end{align}
In the sequel, we will assume that the $\varepsilon_p$ have precisely this expression. For a given $g:\mathcal{Y}\rightarrow\R$, we denote by $N[g]$ the integral 
\begin{align*}
N[g]
  &=\int_{\mathcal{Y}}g(y)N(dy). 	
\end{align*}
With this notation in mind, we have that $\varepsilon_p
  =N[f_p]$, where 
\begin{align}\label{eq:g_pdef}
f_p(q,k)
  &=\Indi{\{p=q\}}k.	
\end{align}

\noindent Additive functionals of $H_n$ can then be expressed as linear combinations of the $\varepsilon_p$, which in turn become functionals of the Poisson process.

\subsection{Mertens' formulas}
\noindent Next we present two classical estimates due to Mertens, which play a central role in the analytic part of our proofs (see~\cite[pp.~16-18]{Tenenbaum2015}). The first Mertens formula provides an approximation for the logarithmically weighted sum of primes
\begin{align}
\Bigg| \sum_{p \in \mathcal{P}_n} \frac{\log p}{p} - \log n - c_0 \Bigg| &\;\le\; \frac{C}{\log n}, \label{eq:mertens1}
\end{align}
while the second gives a refined estimate for the harmonic sum of prime reciprocals
\begin{align}
\Bigg| \sum_{p \le n} \frac{1}{p} - \log\log n - c_1 \Bigg| &\;\le\; \frac{5}{\log n}, \label{eq:mertens2}
\end{align}
for some absolute constants $c_0,c_1$ and $C$. In the harmonic sampling framework, these approximations serve as key tools for evaluating expectations of additive functions and controlling their asymptotic behavior.\\

\subsection{The Dickman distribution}\label{sec:dickmansec}
\noindent 
The Dickman distribution $\mathscr{D}$ was introduced in the seminal work of Dickman~\cite{Dickman1930}, 
arising in the study of the largest prime factor of integers. It is a probability distribution supported on~$[0,\infty)$, whose density is given by $e^{-\gamma}\rho$, where $\gamma$ denotes the Euler-Mascheroni constant and $\rho$ is the classical Dickman function. The function $\rho$ is defined as the unique continuous solution to the delay differential equation
\[
u\,\rho'(u) = -\,\rho(u-1), \quad u > 1, 
\qquad \rho(u) = 1,
\]
for $0 \le u \le 1.$ A probabilistic perspective of this random variable is obtained by regarding the distribution as the unique distribution satisfying the bias identity
\begin{equation} \label{eq:Dickman_bias}
\mathbb{E}[\mathscr{D} f(\mathscr{D})] = \int_0^1 \mathbb{E}[f(\mathscr{D}+t)] \, dt.
\end{equation}
\noindent for all bounded smooth $f : \mathbb{R}^+ \to \mathbb{R}$. The Dickman distribution is infinitely divisible. In particular, we can define a L\'evy process $D=\{D_t;\,t\ge0\}$ such that $D_1$ is equal in law to $\mathscr D$. This process will be referred to as a L\'evy subordinator.\\

\noindent A particular realization of the Dickman distribution can be obtained via a Poisson point process. Let $\eta$ be a Poisson random measure on $\mathcal{X}:=(0,1]^2$ with intensity
\begin{align}\label{eq:muintensity}
\mu(dt,dx)=\frac{\mathbbm{1}_{(0,1]}(x)}{x}\,dt\,dx.	
\end{align}
From here we can define the process
\begin{equation}\label{eq:Dickman_subordinator_from_PPP}
L_t := \int_{\mathcal{X}}\mathbbm{1}_{[0,t]}(s)\,x\,\eta(ds,dx), 
\end{equation}
for $t$ in $[0,1]$. Then $\{L_t;\,t\in[0,1]\}$ is a L\'evy subordinator with L\'evy measure $\mathbbm{1}_{(0,1]}(x)\,dx/x$ and no Gaussian part.
This process will be referred to, in the sequel, as the {Dickman subordinator}.

\medskip\noindent
The above construction naturally induces a wider family of infinitely divisible
random variables obtained as Poisson integrals against~$\eta$.
For a deterministic measurable function $f:(0,1]\to\R$ that is integrable with respect to the measure $x^{-1}dx$ on $(0,1]$, that is,
\[
\int_0^1 \frac{|f(x)|}{x}\,dx<\infty,
\]
we define
\[
X[f] := \int_{\mathcal{X}} f(x)\,\eta(ds,dx).
\]
Then $X[f]$ is infinitely divisible with characteristic function
\begin{align}\label{charfunctionpoissint}
\mathbb{E}\!\left[e^{iuX[f]}\right]
=\exp\!\left\{\int_0^1 \frac{1}{x}\big(e^{iu f(x)}-1\big)\,dx\right\}.	
\end{align}

\section{Main results}\label{sec:main}
\noindent We now present our main results. Recall that any completely additive function $f:\N\rightarrow\R$ is determined by its values in the set $\mathcal{P}$, as the prime decomposition yields the identity 
\begin{align}\label{eq:additivedef}
f(k)
  &=f\left(\prod_{p\in\mathcal{P}}p^{\alpha_p(k)}\right)
  =\sum_{p\in\mathcal{P}}\alpha_p(k)f(p),
\end{align}
for all positive integers $k$. Next we narrow down the type of completely additive functions we will be considering in this paper.\\

\noindent\textit{Completely additive embeddings of real functions.}\\
Relation~\eqref{eq:additivedef} shows that a completely additive arithmetic function is completely determined by its values on primes. In particular, given a real function $\vartheta:\mathbb{R}_{+}\to\mathbb{R}$, we define the induced completely additive function $\iota[\vartheta]:\mathbb{N}\to\mathbb{R}$ by
\[
\iota[\vartheta](p) := \vartheta(\log p),
\]
for all primes $p$. The definition then extends completely additively to all positive integers via the prime factorization rule~\eqref{eq:additivedef}. Equivalently,
\[
\iota[\vartheta](p^\ell)=\ell\,\vartheta(\log p),
\]
for all primes $p$ and non-negative integers $\ell$.\\

\noindent\textit{Specifications of the model}\\
Next we describe the family of choices for $\vartheta$ to which our result applies. Denote by $\mathcal V$ the collection of continuously differentiable functions $\vartheta:\R_+\to\R$ such that $\vartheta$ is eventually non-zero and regularly varying at infinity with positive index. More precisely, for each $\vartheta$ in $\mathcal V$ there exists $\alpha>0$ such that
\[
\frac{\vartheta(cx)}{\vartheta(x)}
\rightarrow c^\alpha
\]
as $x$ tends to infinity, for every $c>0$. We also assume the derivative regularity condition
\[
\frac{1}{|\vartheta(x)|}
\int_1^x \frac{|\vartheta'(s)|}{s}\,ds
\rightarrow 0
\]
as $x$ tends to infinity. Finally, we assume that for every $a>0$ there exist constants $\beta>0$, $C>0$ and $x_0>0$ such that
\[
\left|\frac{\vartheta(ux)}{\vartheta(x)}\right|\le C u^\beta
\]
for all $x\ge x_0$ and all $a/x\leq u\leq 1$. The index $\alpha$ will be denoted by $\partial\vartheta$. In the rest of the section, we will adopt the notation $\mathcal{X}$ and $\eta$, as described in Section \ref{sec:dickmansec}. The main result of the paper is the following.

\begin{theorem}\label{thn:main}
Let $\varphi=\iota[\vartheta]$ for some continuously differentiable $\vartheta$ in $\mathcal V$, and suppose that $\alpha:=\partial\vartheta>0$. Then
\[
\frac{\varphi(H_n)}{\vartheta(\log n)}
\;\;\xrightarrow{\;\mathrm{Law}\;}\;\;
\mathcal L(
\int_{\mathcal X}x^\alpha\eta(ds,dx)
 |\;
\int_{\mathcal X}x\eta(ds,dx)\leq 1),
\]
as $n$ tends to infinity.
\end{theorem}

\noindent The next proposition is the main ingredient for the proof of \eqref{thn:main}.  Its content is of independent interest. Let $\varepsilon_p$ denote the random variables introduced in Section \ref{eq:epsilonpdef}. 

\begin{proposition}\label{prop:main}
Let $\varphi$ be a completely additive function of the form $\varphi = \iota[\vartheta]$ for some continuously differentiable $\vartheta $ in $ \mathcal{V}$, and suppose that $\alpha:=\partial\vartheta>0$. For every $u_1,u_2$ in $\mathbb R$, consider the random variable
\[
Z_n \;=\; \sum_{p \in \mathcal{P}_n} \left( \frac{u_1 \,\varphi(p)}{\vartheta\circ\log(n)} \;+\; u_2 \,\frac{\log p}{\log n} \right)\, \varepsilon_p.
\]
Then, as $n$ tends to infinity,
\[
Z_n \;\;\xrightarrow{\;\mathrm{Law}\;}\;\;
u_1  \int_{\mathcal{X}}x^{\alpha}\eta(ds,dx) + u_2 \int_{\mathcal{X}}x\eta(ds,dx).
\]
\end{proposition}

As anticipated, this proposition serves as the key building block for determining the behavior of the additive functionals discussed in Section~\ref{sec:intro}.\\

\noindent To the best of our knowledge, this type of Dickman-conditioned Poisson integral limit for completely additive functions under harmonic sampling has not been previously formulated in this form. As mentioned in the introduction, most previous work concerns the case where $H_n$ is replaced by a uniform sample from $[n]$. In that setting, Billingsley's results yield infinitely divisible limits of a different nature from the Dickman-type conditional limit exhibited in Theorem~\ref{thn:main}.

\section{Proof of Proposition \ref{prop:main}}
\noindent We start representing the distribution of the variables \eqref{eq:epsilondef} in terms of the point process $N$. To this end, we consider the representation $\varepsilon_p=N[f_p]$, with $f_p$ given as in \eqref{eq:g_pdef}.\\

\noindent Let $\lambda$ in $\R$ be given and define $v_i=\lambda u_i$. In order to determine the asymptotic behavior of the law of $Z_n$, it suffices to determine the asymptotics of 
\[
E\left[
\exp\left\{ \mathbf{i} \sum_{p \in\mathcal{P}_n} \left( \frac{v_1 \,\varphi(p)}{\vartheta \circ \log(n)} +  \frac{v_2\log p}{\log n} \right) \varepsilon_p \right\}
\right].
\]
To this end, we will make use of the observation at the end of Section \ref{sec:dickmansec}, which states that Poisson integrals have an explicit characteristic function that can be written in terms of the underlying control measure. To carry out this program, we observe that with the above notation in mind, 
\[
\lambda Z_n
  = \int_{\mathcal{Y}} g_n(z)N(dz), 
\]
where 
\[
g_n(z)
  = \sum_{p \in\mathcal{P}_n} \left( \frac{v_1 \,\varphi(p)}{\vartheta\circ \log(n)} +  \frac{v_2\log p}{\log n} \right) f_p(z).
\]

\noindent Then
\[
E\big[\exp(i\lambda Z_n)\big] = \exp\left\{ \int_{\mathcal{Y}} \big(e^{\mathbf{i} g_n(z)} - 1\big) \, \nu(dz) \right\}.
\]
where $\nu$ is given as in Section \ref{eq:epsilonpdef}. We are then left with the problem of determining the asymptotic behavior of
\[
I_n := \int_{\mathcal{Y}} \big(e^{\mathbf{i} g_n(z)} - 1\big) \, \nu(dz).
\]
Observe that if $z=(q,k)$, then
\[
g_n(q,k)
=
\sum_{p\in\mathcal P_n}
\left(
\frac{v_1\varphi(p)}{\vartheta\circ\log(n)}
+
\frac{v_2\log p}{\log n}
\right)f_p(q,k)
=
k\mathbbm 1_{\{q\le n\}}
\left(
\frac{v_1\varphi(q)}{\vartheta\circ\log(n)}
+
\frac{v_2\log q}{\log n}
\right).
\]
From here it follows that
\[
I_n 
=
\sum_{k \geq 1} \sum_{q \in \mathcal{P}_n}
\frac{1}{kq^k}
\left(
e^{\mathbf{i} k \left( \frac{v_1 \varphi(q)}{\vartheta\circ \log(n)} + \frac{v_2 \log q}{\log n} \right)} - 1
\right).
\]
Next we show that the term $k$ equal to one yields the main contribution as $n$ tends to infinity. To this end, we consider the decomposition
\[
I_n = \tilde{I}_n + R_n,
\]
where 
\[
\tilde{I}_n
  := \sum_{q \in \mathcal{P}_n} \left( e^{ i \left( \frac{ v_1\varphi(q)}{\vartheta\circ \log(n)} + \frac{v_2 \log(q)}{\log(n)} \right)} - 1 \right) \frac{1}{q}
\]
and 
\[
R_n
  := \sum_{k \geq 2} \sum_{q \in \mathcal{P}_n} \left( e^{ i k \left( \frac{v_1 \varphi(q)}{\vartheta\circ \log(n)} + \frac{v_2 \log(q)}{\log(n)} \right)} - 1 \right) \frac{1}{k q^k}
\]
For each fixed $q$ and $k\ge2$,
\[
k\left(
\frac{v_1\varphi(q)}{\vartheta\circ\log(n)}
+
\frac{v_2\log q}{\log n}
\right)\to0.
\]
Moreover,
\[
\left|
\left(
e^{ i k \left( \frac{v_1 \varphi(q)}{\vartheta\circ \log(n)} + \frac{v_2 \log q}{\log n} \right)} - 1
\right)
\frac{1}{kq^k}
\right|
\le \frac{2}{kq^k},
\]
and
\[
\sum_{k\ge2}\sum_{q\in\mathcal P}\frac{1}{kq^k}<\infty.
\]
Thus, by dominated convergence, $R_n\to0$. It then suffices to handle $\tilde I_n$.  To this end, we observe that 
\[
\tilde{I}_n
  = \int_{2}^n \left( e^{i\left(\frac{v_1\varphi(t)}{\vartheta\circ \log(n)} + \frac{v_2 \log(t)}{\log(n)}\right)} - 1\right) F(dt)
\]
where $F$ is the cumulative prime-reciprocal function
\[
F(x) = \sum_{q \in \mathcal{P}_{\floor{x}}} \frac{1}{q}.
\]
The fact that $\varphi=\iota[\vartheta]$ allows us to replace $\varphi(t)$ by $\vartheta\circ\log(t)$ at the atoms of $F$. By Mertens' approximation \eqref{eq:mertens2}, if
\[
G(x):=F(x)-\log\log x-c_1,
\]
then $|G(x)|\le C/\log x$ for $x\ge2$. We therefore consider the decomposition
\[
\tilde{I}_n = \hat{I}_n + \hat{R}_n,
\]
where
\begin{align*}
\hat{I}_n
  &= \int_{2-}^n \left( e^{ i \left( \frac{v_1 \vartheta\circ \log(t)}{\vartheta\circ \log(n)} + \frac{v_2 \log(t)}{\log(n)} \right) } - 1 \right) \log \log (dt)
\\
\hat{R}_n
  &= \int_{2-}^n \left( e^{ i \left( \frac{v_1 \vartheta\circ \log(t)}{\vartheta\circ \log(n)} + \frac{v_2 \log(t)}{\log(n)} \right) } - 1 \right) G(dt).	
\end{align*}
The treatment of the terms $\hat{I}_{n}$ and $\hat{R}_n$ will be carried in separate steps, which will be followed by an identification of the corresponding limit.\\

\noindent\textbf{Step I}\\
We begin handling the term $\hat{R}_n$. To this end, define
\[
\phi(t):=\frac{v_1\vartheta\circ \log(t)}{\vartheta\circ \log(n)}+\frac{v_2\log t}{\log n},
\qquad 
f(t):=e^{\mathbf{i}\phi(t)}-1 ,
\]
so that
\[
\widehat R_n
  = \int_{(2,n]} f(t)\,G(dt).
\]
By integration by parts,
\[
\int_{(2,n]} f(t)\,G(dt)
=
\bigl[f(t)G(t)\bigr]_{2-}^{\,n}
-\int_{2}^{n} G(t)\,df(t).
\]
Since \(df(t)=\mathbf{i}e^{\mathbf{i}\phi(t)}\,\phi'(t)\,dt\), we obtain
\[
\widehat R_n
=
\left( e^{\mathbf{i}\phi(n)}-1 \right)G(n)
-\left( e^{\mathbf{i}\phi(2-)}-1 \right)G(2-)
-\mathbf{i}\int_{2}^{n} G(t) e^{\mathbf{i}\phi(t)}\,\phi'(t)\,dt.
\]
Using $|e^{ix}-1|\le 2\wedge |x|$ and Mertens' approximation \eqref{eq:mertens2}, we have $|G(t)|\le C/\log t$ for $t\ge2$. Hence
\[
|\widehat R_n|
\le
C|f(n)|/\log n
+
C|f(2-)|
+
C\int_2^n \frac{|\phi'(t)|}{\log t}\,dt .
\]
The first term tends to zero since $|f(n)|\le2$. Moreover,
\[
|f(2-)|\le C|\phi(2)|
\le
C\left(\frac{1}{|\vartheta\circ\log(n)|}+\frac{1}{\log n}\right),
\]
which also tends to zero. Finally,
\[
\int_2^n \frac{|\phi'(t)|}{\log t}\,dt
\le
\frac{C}{|\vartheta\circ\log(n)|}
\int_2^n \frac{|(\vartheta\circ\log)'(t)|}{\log t}\,dt
+
\frac{C}{\log n}\int_2^n \frac{dt}{t\log t}.
\]
The second term on the right-hand side is $O(\log\log n/\log n)$ and therefore tends to zero. For the first term, the change of variables $s=\log t$ gives
\[
\frac{1}{|\vartheta(\log n)|}
\int_{\log 2}^{\log n}\frac{|\vartheta'(s)|}{s}\,ds,
\]
which tends to zero by the regular-variation assumption on $\vartheta$ together with the derivative regularity imposed in the definition of $\mathcal V$. Consequently, $\widehat R_n\to0$.\\

\noindent\textbf{Step II}\\
It now remains to determine the asymptotic behavior of $\hat{I}_n$. To this end, we observe that
\[
\hat{I}_n = \int_{2}^{n} \left( \exp\left\{ i \left( \frac{v_1 \vartheta\circ \log(t)}{\vartheta\circ \log(n)} + \frac{v_2 \log(t)}{\log(n)} \right) \right\} - 1 \right) \frac{1}{t \log(t)} dt.
\]
By the change of variables $u=\log(t)/\log(n)$, we obtain
\begin{align*}
\hat{I}_n
  &=\int_{\log(2)/\log(n)}^1
  \left(
  e^{ i \left( \frac{ v_1 \vartheta (u\log n) }{ \vartheta\circ \log(n) } + v_2 u \right)  } - 1
  \right)\frac{du}{u}.
\end{align*}
Since $\vartheta$ is regularly varying with index $\alpha=\partial\vartheta>0$, we have
\[
\frac{\vartheta(u\log n)}{\vartheta(\log n)}\rightarrow u^\alpha
\]
for every $u$ in $(0,1]$, and the convergence is uniform on compact subsets of $(0,1]$. We split the integral over $(\log(2)/\log(n),\delta)$ and $[\delta,1]$. On $[\delta,1]$, the uniform convergence theorem for regularly varying functions gives
\[
\int_{\delta}^{1}
\left(
e^{ i \left( \frac{ v_1 \vartheta (u\log n) }{ \vartheta\circ \log(n) } + v_2 u \right)  } - 1
\right)\frac{du}{u}
\longrightarrow
\int_{\delta}^{1}
\left(
e^{i(v_1u^\alpha+v_2u)}-1
\right)\frac{du}{u}.
\]
On $(\log(2)/\log(n),\delta)$, the domination condition in the definition of $\mathcal V$ yields, for some $\beta>0$,
\[
\left|\frac{\vartheta(u\log n)}{\vartheta(\log n)}\right|\le C u^\beta.
\]
Thus, using $|e^{ix}-1|\le |x|$,
\[
\left|
e^{ i \left( \frac{ v_1 \vartheta (u\log n) }{ \vartheta\circ \log(n) } + v_2 u \right)  } - 1
\right|
\le C(u^\beta+u).
\]
Since $(u^\beta+u)/u$ is integrable near zero, the contribution of the interval $(\log(2)/\log(n),\delta)$ is bounded uniformly by a quantity that tends to zero as $\delta$ tends to zero. Therefore
\[
\lim_n I_n=\lim_n \hat{I}_n = \int_0^1 \left( \exp\{i(v_1 u^\alpha + v_2u)\} - 1 \right) \frac{du}{u}.
\]

\noindent \textbf{Step III}\\
Let $\mathcal{X}$, $\mu$ and $\eta$ be the space, control measure and Poisson point process introduced in Section \ref{sec:dickmansec}. By the characteristic function formula for Poisson integrals, the random variable
\[
v_1\int_{\mathcal X}x^\alpha\eta(ds,dx)+v_2\int_{\mathcal X}x\eta(ds,dx)
\]
has characteristic function
\[
\exp\left\{\int_0^1
\left(
e^{\mathbf{i}(v_1x^\alpha+v_2x)}-1
\right)\frac{dx}{x}
\right\}.
\]
The limit obtained in Step II is precisely this exponent. Therefore
\[
\lim_{n\to\infty}\E[e^{\mathbf{i}\lambda Z_n}]
=
\E\left[
\exp\left\{
\mathbf{i}\lambda u_1\int_{\mathcal X}x^\alpha\eta(ds,dx)
+
\mathbf{i}\lambda u_2\int_{\mathcal X}x\eta(ds,dx)
\right\}
\right],
\]
and the proof of Proposition \ref{prop:main} follows.

\section{Proof of Theorem \ref{thn:main}}
\noindent Let $\Psi:\R\rightarrow\R$ be a smooth bounded test function. To prove the result, it suffices to analyze the asymptotic behavior of the expectation of $\Psi(\varphi(H_n)/\vartheta\circ\log(n))$. In view of the identity \eqref{eq:HN_decomposition}, we have that 
\begin{equation*}
\E\left[\Psi\left(\frac{\varphi(H_n)}{\vartheta\circ\log(n)}\right)\right] =
\E\Big[\Psi\big(\sum_{p \in \mathcal{P}_{n}} \varphi(p^{\varepsilon_{p}})/\vartheta\circ\log(n)\big)\Indi{A_{n}}\Big]/\Pb[A_n].
\end{equation*}
Next, by writing $A_n=\left\{Z_{n,2}\leq 1\right\}$, and using the complete additivity identity $\varphi(p^{\varepsilon_p})=\varepsilon_p\varphi(p)$, with
\begin{align*}
Z_{n,2}
  &:=\frac{1}{\log(n)}\sum_{p\in\mathcal{P}_n}\log(p)\varepsilon_p,
\end{align*}
we get that 
\begin{align*}
\E\left[\Psi\left(\frac{\varphi(H_n)}{\vartheta\circ\log(n)}\right)\right]
  &=\E[\tilde{\Psi}(Z_{n,1},Z_{n,2})]/\E[\Indi{\{Z_{n,2}\leq 1\}}],	
\end{align*}
with $\tilde{\Psi}(a,b):=\Psi(a)\Indi{\{b\le 1\}}$ and 
\begin{align*}
Z_{n,1}
  &:=\frac{1}{\vartheta\circ\log(n)}\sum_{p \in \mathcal{P}_{n}} \varphi(p)\varepsilon_p.
\end{align*}
By Proposition \ref{prop:main}, the vector $(Z_{n,1},Z_{n,2})$ converges in law, in the sense of arbitrary linear combinations, to 
\begin{align*}
\left(\int_{\mathcal{X}}x^{\alpha}\eta(ds,dx), \int_{\mathcal{X}}x\eta(ds,dx)\right).
\end{align*}
Since the second limiting coordinate has the Dickman distribution, it has no atom at $1$. Therefore
\[
\Pb\left[\int_{\mathcal X}x\eta(ds,dx)=1\right]=0.
\]
Hence the set $\{(a,b):b\leq 1\}$ is a continuity set for the limiting vector, and the result follows from the convergence in law of $(Z_{n,1},Z_{n,2})$.\\

\noindent \textbf{Acknowledgements}\\
Arturo Jaramillo Gil was supported by
the grant CBF2023-2024-2088.

\bibliographystyle{plain}
\bibliography{references.bib}

@book{Tenenbaum2015,
  author    = {G{\'e}rald Tenenbaum},
  title     = {Introduction to Analytic and Probabilistic Number Theory},
  series    = {Cambridge Studies in Advanced Mathematics},
  volume    = {163},
  publisher = {Cambridge University Press},
  address   = {Cambridge},
  year      = {2015},
  edition   = {Third edition},
  note      = {Translated from the 4th French edition by Patrick D. F. Ion},
  isbn      = {978-1-107-08964-2},
  doi       = {10.1017/CBO9781107415801}
}

@article{RenyiTuran1958,
    AUTHOR = {R\'enyi, A. and Tur\'an, P.},
     TITLE = {On a theorem of {E}rd\"os-{K}ac},
   JOURNAL = {Acta Arith.},
  FJOURNAL = {Polska Akademia Nauk. Instytut Matematyczny. Acta Arithmetica},
    VOLUME = {4},
      YEAR = {1958},
     PAGES = {71--84},
      ISSN = {0065-1036},
   MRCLASS = {10.00},
  MRNUMBER = {96629},
MRREVIEWER = {P.\ Erd\H os},
       DOI = {10.4064/aa-4-1-71-84},
       URL = {https://doi.org/10.4064/aa-4-1-71-84},
}

@article{Dickman1930,
  author  = {Dickman, Karl},
  title   = {On the frequency of numbers containing prime factors of a certain relative magnitude},
  journal = {Arkiv f{\"o}r Matematik, Astronomi och Fysik},
  volume  = {22A},
  number  = {10},
  pages   = {1--14},
  year    = {1930}
}

@article{JaramilloYang2023,
  author  = {Arturo Jaramillo and Xiaochuan Yang},
  title   = {Approximation of Smooth Numbers for Harmonic Samples: A Stein Method Approach},
  journal = {Preprint},
  year    = {2023},
  note    = {arXiv:2307.00000}
}

@misc{BhattacharjeeSchulte2022,
  author       = {Chinmoy Bhattacharjee and Matthias Schulte},
  title        = {Dickman Approximation of Weighted Random Sums in the Kolmogorov Distance},
  year         = {2022},
  archivePrefix= {arXiv},
  eprint       = {2211.10171},
  primaryClass = {math.PR},
  url          = {https://arxiv.org/abs/2211.10171}
}

@book{Elliott1979,
    AUTHOR = {Elliott, P. D. T. A.},
     TITLE = {Probabilistic number theory. {I}},
    SERIES = {Grundlehren der Mathematischen Wissenschaften},
    VOLUME = {239},
      NOTE = {Mean-value theorems},
 PUBLISHER = {Springer-Verlag, New York-Berlin},
      YEAR = {1979},
     PAGES = {xxii+359+xxxiii pp. (2 plates)},
      ISBN = {0-387-90437-9},
   MRCLASS = {10-02 (10K20 60C05)},
  MRNUMBER = {551361},
MRREVIEWER = {J.\ Kubilius},
}

@book{Elliott1980,
  author    = {Peter D. T. A. Elliott},
  title     = {Probabilistic Number Theory II: Central Limit Theorems},
  series    = {Grundlehren der Mathematischen Wissenschaften},
  volume    = {240},
  publisher = {Springer-Verlag},
  address   = {Berlin--New York},
  year      = {1980}
}

@article{KowalskiNikeghbali2014,
  author  = {E. Kowalski and A. Nikeghbali},
  title   = {Mod-discrete expansions in probability and number theory},
  journal = {Probability Theory and Related Fields},
  volume  = {158},
  number  = {3--4},
  year    = {2014},
  pages   = {859--885}
}

@article{ChenJaramilloYang2021,
  author  = {L.~H.~Y. Chen and A. Jaramillo and X.~C. Yang},
  title   = {A probabilistic approach to the Erd{\H{o}}s--Kac theorem for additive functions},
  journal = {Preprint (2021)},
  note    = {arXiv:2102.05094}
}

@book{Kubilius1964,
  author    = {A. Kubilius},
  title     = {Probabilistic Methods in the Theory of Numbers},
  publisher = {American Mathematical Society},
  series    = {Translations of Mathematical Monographs},
  volume    = {11},
  year      = {1964}
}

@article{ErdosKac1940,
  author  = {P. Erd{\H{o}}s and M. Kac},
  title   = {The Gaussian law of errors in the theory of additive number theoretic functions},
  journal = {American Journal of Mathematics},
  volume  = {62},
  year    = {1940},
  pages   = {738--742}
}

@article{Billingsley1974,
  author  = {P. Billingsley},
  title   = {The probability theory of additive arithmetic functions},
  journal = {Annals of Probability},
  volume  = {2},
  year    = {1974},
  pages   = {749--791}
}

\end{document}